\documentclass[reqno,12pt]{amsart}
\usepackage{bbm}
\usepackage{amsmath}
\usepackage{txfonts}
\usepackage{stmaryrd}
\usepackage{amssymb}
\usepackage{amsfonts}
\usepackage{times}
\usepackage{mathrsfs}

\def\D{\mathbb{D}}

\def\C{\mathbb{C}}

\newcommand{\DD}{{\mathbb D}}

\def\msk{\medskip}

\def\bege{\begin{equation}} \def\ende{\end{equation}}
\def\begr{\begin{eqnarray}} \def\endr{\end{eqnarray}}

\def\square{\vbox{
    \hrule height .4pt
    \hbox{\vrule width .4pt height 7pt \kern 7pt
       \vrule width .4pt}
    \hrule height .4pt }}

\date{}

\begin{document}

\title{On the structure of invariant subspaces\\ for the shift operator on Bergman spaces}

\author{Junfeng Liu}

\address{Junfeng Liu\\ Faculty of Information Technology, Macau University of Science and Technology, Avenida Wai Long, Taipa, Macau. }   \email{jfliu997@sina.com}

\thanks{This project was partially supported by the Macao Science and Technology Development Fund (No.186/2017/A3).}

\maketitle
\begin{abstract}It is well known that the structure of nontrival invariant subspaces for the shift operator
on the Bergman space is extremely complicated and very little is known about their specific structures,
and that a complete description of the structure seems unlikely.
In this paper, we find that any invariant subspace $M(\neq \{0\})$
for the shift operator $M_z$ on the Bergman space $A^{^p}_\alpha(\D)~(1\leq p< \infty,~~-1<\alpha<\infty)$ contains a nonempty subset that lies in
$A^{^p}_\alpha(\D)\backslash H^{^p}(\D)$. To a certain extent, this result describes the specific
structure of every nonzero invariant subspace for the shift operator $M_z$ on $A^{^p}_\alpha(\D)$.
\end{abstract}



\section{Introduction and Preliminaries}

One of the most famous unsolved problems in mathematics is the invariant subspace problem (cf. \cite{pda04}, p.268).
In the Hardy space $H^2(\D)$, A. Beurling \cite{ab49} found that every invariant subspace for the shift operator
has an very elegant description. Also, his theory extends with litter change to the Hardy space $H^{^p}(\D)~(1\leq p< \infty).$

In order to state Beurling's theorem for the Hardy space, we first recall
the notion of the inner function.
If $\phi(z)$ is a functional in $H^\infty(\D)$ such that
$|\widetilde{\phi}(e^{i\theta})|=1$ a.e. on $\partial\D$,
where $\widetilde{\phi}(e^{i\theta})=\lim_{r\rightarrow 1^-}\phi(re^{i\theta})$
a.e. on $\partial\DD$,
then $\phi$ is called an inner function.

\msk
\noindent
\textbf{Theorem A}
\textbf{(Beurling's theorem).}
\emph{Every invariant subspace $M(\neq \{0\})$ for
the shift operator $M_z$ on the Hardy space $H^{^p}(\D)\ ~~(1\leq p<\infty)$ has the form
$$M=\phi H^{p}(\D)$$
for some inner function $\phi$
(for more details, see \cite{ab49}, \cite{pd70}, \cite{pda04}).}
\msk

It is well known that in the first half of twentieth century,
the research interest of complex analysis dipped to a
low point. It was Beurling's work \cite{ab49} that revived the research interest
of complex analysis, and that opened a new research direction of the invariant subspace problem.

It is worth noticing that Beurling's theorem explicitly describe
the specific structure of every nontrivial invariant subspace for the shift operator
$M_z$ on the Hardy space $H^{^p}(\DD)~~(1\leq p< \infty)$,
but in contrast it is still not known how to formulate the specific structure of nontrivial invariant subspaces for the shift operator on the
Bergman space $A^{^p}_\alpha(\D)~~(1\leq p<\infty,~~-1<\alpha<\infty)$.

In the topic of the invariant subspace for the shift operator on the
Bergman space, the most famous result is due to A. Aleman, S. Richter and C. Sundberg \cite{asc96} so far.
Later  S. M. Shimorin \cite{sms01} and \cite{sms02}
found simpler proofs and extended the
result of \cite{asc96} from $A^2(\D)$ to $A^2_\alpha(\D)~~(-1<\alpha\leq 1)$,
in which the extension of $-1<\alpha< 0$ is due to \cite{sms01}, while the extension of
$0<\alpha\leq 1$ is due to \cite{sms02}.
Indeed, S. M. Shimorin \cite{sms01} and \cite{sms02} obtained  more general theorems for Hilbert spaces.
Now, we state these results.

\msk
\noindent
\textbf{Theorem B}
\textbf{(the Aleman-Richter-Sundberg theorem \cite{asc96}).}
\emph{Every invariant subspace $M$ for the shift operator $M_z$ on the Bergman space $A^2_\alpha(\D)~(-1<\alpha\leq1)$
has the form
$$M=[M\ominus zM],$$
where $[M\ominus zM]$ denotes the smallest invariant subspace for $M_z$
that contains $M\ominus zM$
(for more details, see \cite{asc96}, \cite{pda04}, \cite{hkz2000}, \cite{sms01} and \cite{sms02}).}
\msk

This is an interesting result and can be regarded as a version of Beurling's theorem for the Bergman space.
On the other hand, it follows from the existence result in \cite{abf85} and explicit examples in \cite{hh93} and \cite{hrs96} that the
dimension of the linear subspace $M\ominus zM$ of $A^2_\alpha(\D)$ in Theorem B may be any positive integer or the infinity (see also \cite{pda04}
and \cite{hkz2000}).
Moreover, it is not even known how or whether the result in Theorem B
might extend to $A^p_\alpha(\D)$ for other values of $p$ and $\alpha$
(cf. \cite{pda04} p.3, \cite{hkz2000} p.187, \cite{sms01} and \cite{sms02}).

\msk

But for
the nontrivial invariant subspace $M$ for the shift operator $M_z$
on the Bergman space, very little
is known about their specific structures (cf. \cite{kzb07} p.95).
In other words, the specific structure of nontrivial invariant subspaces
for the shift operator on the Bergman space has not been
characterized (cf. \cite{pda04} p.3).

For a long time, it has been known that the structure of nontrival invariant subspaces for the shift operator
on the Bergman space is extremely complicated,
and that a complete description of the structure seems unlikely
(cf. \cite{asc96} p.277 and \cite{pda04} p.3, see also \cite{abf85}, \cite{hh93} and \cite{hrs96}).
In any case, the most intriguing topic in the theory of Bergman spaces is the
structure problem of invariant subspaces for the shift operator
(cf. \cite{kzb07} p.95).

In this paper, we find that any invariant subspace $M(\neq \{0\})$
for the shift operator $M_z$ on the Bergman space $A^{^p}_\alpha(\D)~(1\leq p< \infty,~~-1<\alpha<\infty)$ contains a nonempty subset that lies in
$A^{^p}_\alpha(\D)\backslash H^{^p}(\D)$. To a certain extent, this result characterizes the specific
structure of every invariant subspace for the shift operator $M_z$ on $A^{^p}_\alpha(\D)$.

To this end, we recall the basic concept of the Bergman space from  \cite{pda04}, \cite{hkz2000},
\cite{kzb07}, which are used in the main results.

Let $\D$ and $\partial\D$ denote the open unit disk and the unit circle on the complex plane respectively.
The (weighted) Bergman space $A^{^p}_\alpha(\D)~(1\leq p< \infty,~-1<\alpha<\infty)$ is defined by

$$A^{^p}_\alpha(\D)=\Big\{f:\mbox{the~ function~} f(z)~ {\mbox{is analytic in}}~~ \D~~ {\mbox{and}} ~~\int_{\D} |f(z)|^{^p} dA_{\alpha}(z)< \infty\Big\}$$
with the norm

$$\|f\|_{_{A^{^p}_\alpha(\D)}}=\left(\int_{\D} |f(z)|^{^p} dA_{\alpha}(z)\right)^\frac{1}{p},~~f\in A^{^p}_\alpha(\D),$$
where $dA_{\alpha}(z)= \frac{\alpha+1}{\pi}(1-|z|^2)^{\alpha}dA(z)$, while $dA$
denotes the planar Lebesgue measure (the area measure) on $\D$ with $dA(\D)=\pi$.
When $\alpha=0$, the  (weighted) Bergman space $A^{^p}_0(\D)~~(1\leq p< \infty)$ is just the (classical) Bergman space $A^{^p}(\D)$.
Sometimes $A^{^p}_0(\D)$ is written as $L^{p}_{a}(dA\alpha)$.
\msk

Now, we gather together some notions of the shift operator and the invariant subspace for shift operators on Hardy spaces
and Bergman spaces.

\msk
\noindent
\textbf{Definition 1.}
\emph{Let $X$ be a Banach space of analytic
functions in $\DD$ (for example, the Hardy space $H^{^p}(\DD)~~(1\leq p< \infty)$
and the Bergman space $A^{^p}_\alpha(\D)~( 1\leq p<\infty,~~-1<\alpha<\infty)$).
Let $M$ be a closed linear subspace of $X$.}

\emph{(1). If the operator $M_z$ on X is defined by
$$(M_z)(z)=zf(z),~~~~f\in X ,$$
then $M_z$ is called the shift operator on $X$.}

\emph{(2). If $(M_zM=)\ zM \subset M$, then $M$ is called an invariant subspace for the shift operator
$M_z$ on $X$,
or is briefly called an invariant subspace in $X$.}

\msk

\section{Main Results}

It is well known that as linear spaces, $H^{^p}(\D)~(1\leq p<\infty)$ is a linear subspace of
$A^{^p}_\alpha(\D)~(1\leq p< \infty,~~-1<\alpha<\infty)$ and the inverse is not true, that is,
$H^{^p}(\D)\subsetneqq A^{^p}_\alpha(\D)$, which can be seen from
\cite{pda04} p.74, p.77, p.80, and \cite{kzb07} p.100 Problem 28 as well as 30 and so on.
But as normed space, the norms on $H^{^p}(\D)$ and $A^{^p}_\alpha(\D)$ are different,
and what is more, now there is a big difference between the structure theory of invariant subspaces for the shift
operator $M_z$ on $H^{^p}(\D)$ and $A^{^p}_\alpha(\D)$.
First of all, by Beurling's theorem for the Hardy space, the specific structure of every invariant subspace
for the shift operator $M_z$ on  $H^{^p}(\D)~(1\leq p<\infty)$ is very clear.
On the other hand, as mentioned in the introduction, very little is
known about the specific structure of nontrivial invariant subspaces
for the shift operator $M_z$ on $A^{^p}_\alpha(\D)$.

To a certain extent, the following result describes the specific
structure of every invariant subspace for the shift operator $M_z$ on $A^{^p}_\alpha(\D)$.

\msk

\noindent\textbf{
Theorem 1.}
\emph{Let $M(\neq 0)$ be any invariant
subspace for the shift operator $M_z$ on the Bergman space
$A^{^p}_\alpha(\D)~~(1\leq p< \infty,~~-1<\alpha<\infty)$. If $M$ is  contained in $H^{^p}(\D)~~(1\leq p< \infty)$,
then $M$
has the form
$$M=\phi H^{^p}(\D)$$
for some inner function $\phi$.}

\emph{\textbf{Proof.}}
To make a difference, the shift operator on the Hardy space $H^{^p}(\D)$ is now written as $S$.
Then for any $f\in H^{^p}(\D)$, it is clear that $Sf=M_z f$, where $M_z$ denotes the shift operator
on the Bergman space
$A^{^p}_\alpha(\D)$.
Since $M\subset H^{^p}(\D) \subset A^{^p}_\alpha(\D)$, it follows that
\begin{equation}
  SM=M_z M\subset M.
\end{equation}

We now show that $M$ is a closed set in the Hardy space $H^{^p}(\D)$. In fact, let $\{f_n\}$
be any sequence in $M$ and $f\in H^{^p}(\D)$ such that
$$\|f_n-f\|_{_{H^{^p}(\D)}}\rightarrow 0$$
as $n\rightarrow \infty$.
Since $f_n\in M\subset H^{^p}(\D) \subset A^{^p}_\alpha(\D)$ and $f\in H^{^p}(\D) \subset A^{^p}_\alpha(\D)$,
it can be obtained that

$$\|f_n-f\|_{_{A^{^p}_\alpha (\D)}}
\leq \|f_n-f\|_{_{H^{^p}(\D)}}\rightarrow 0$$
as $n\rightarrow \infty.$
This implies the sequence $\{f_n\}$ in $M$ converges to $f$ in the norm of $A^{^p}_\alpha(\D)$.
Since $M$ is a closed set in the Bergman space $A^P_\alpha(\D)$, it follows that $f\in M$. It can be seen from
the above result that $M$ is a closed set in the Hardy space $H^{^p}(\D)$. Therefore by (1), $M(\neq \{0\})$
is a invariant subspace for the shift operator $S$ on the Hardy space
$H^{^p}(\D)$. Thus by Beurling's theorem for the Hardy space, $M$ has the form
$$M=\phi H^{^p}(\D)$$
for some inner function $\phi$.
This completes the proof.

\msk

\noindent\textbf{
Theorem 2.}
\emph{Any invariant subspace $M(\neq \{0\})$ for the shift operator $M_z$ on the Bergman space $A_\alpha^{^p}(\D)~(1\leq p<\infty, -1<\alpha <\infty)$
contains a nonempty subset that lies in $A_\alpha^{^p}(\D)\backslash H^{^p}(\D)$.}

\emph{\textbf{Proof.}}
Since $H^{^p}(\D) \subset A_\alpha^{^p}(\D)$ for all $1\leq p<\infty$ and all $-1<\alpha <\infty$, it follows
that if the conclusion of this theorem were not true, then there would be an invariant subspace $M(\neq \{0\})$ for the shift operator $M_z$
on some $A^{^p}_\alpha(\D)$ such that $M\subset H^{^p}(\D)$. Thus by Theorem 1, $M$ has the from
$$M=\phi H^{^p}(\D)$$
for some inner function $\phi$. Therefore $\phi H^{^p}(\D)(=M)$ would be a closed set in $A^{^p}_\alpha(\D)$ (in the norm of $A^{^p}_\alpha(\D)$).

we now show that for any inner function $\phi$, it is true that $\phi H^{^p}(\D)$ is not a closed set in $A^{^p}_\alpha(\D)$ (in the norm of $A^{^p}_\alpha(\D)$).
In fact, it is well known that $H^{^p}(\D) \subsetneqq A^{^p}_\alpha(\D)$.
Let $f_{_{0}}$ be an arbitrary function in $A^{^p}_\alpha(\D)\backslash H^{^p}(\D)$.
Since the set of all polynomials is dense in $A^{^p}_\alpha(\D)$ (in the norm of $A^{^p}_\alpha(\D)$), there is a sequence $\{Q_n\}$
of polynomials such that
$$\|Q_n - f_{_{0}}\|_{_{A^{^p}_\alpha(\D)}}\rightarrow 0$$
as $n\rightarrow \infty$.
From this we can obtain
\begin{eqnarray*}
  \|\phi Q_n- \phi f_{_{0}}\|^p_{A^p_\alpha(\D)}&=&\frac{\alpha+1}{\pi}\int_{\D} |\phi(z)Q_n(z)-\phi(z)f_{_{0}}(z)|^p(1-|z|^2)^\alpha dA(z) \\
   &\leq& \frac{\alpha+1}{\pi}\|\phi\|^p_{_{H^\infty(\D)}}\int_{\D} |Q_n(z)-f_{_{0}}(z)|^p (1-|z|^2)^\alpha dA(z) \\
   &=&\|\phi\|^p_{_{H^\infty(\D)}}\|Q_n-f_{_{0}}\|^p_{_{A^p_\alpha(\D)}} \rightarrow 0
\end{eqnarray*}
as $n\rightarrow \infty$. Since $\phi \in H^{\infty}(\D)$ and $f_0 \in A^{^p}_\alpha(\D)$,
it follows that $\phi f_0\in A^{^p}_\alpha(\D)$. Moreover, it
 is clear that $\phi Q_n\in \phi H^{^p}(\D)$. Thus it remains to show that
$\phi f_{_{0}}\notin \phi H^{^p}(\D)$. Indeed, since $\phi\in H^{\infty}(\D)$, it follows that $\phi H^{^p}(\D)\subset H^{^p}(\D)$. Thus
 if $\phi f_{_{0}} \in \phi H^{^p}(\D)$, then $\phi f_{_0} \in H^{^p}(\D)$.

Let $g(z)=\phi(z)f_{_0}(z)$ for all $z\in\D$. Then $g\in H^p(\D)$,
and
$$\lim_{r\rightarrow 1^-}f_{_0}(re^{i\theta})=\lim_{r\rightarrow 1^-}
\frac{g(re^{i\theta})}{\phi(re^{i\theta})}=\frac{\widetilde{g}(e^{i\theta})}{\widetilde{\phi}(e^{i\theta})}~~~ \mbox{a.~e. ~on} ~~~ \partial\D,$$
where $\widetilde{g}(e^{i\theta})=\lim_{r\rightarrow1^-}g(re^{i\theta})$ a.e. on $\partial\D$ and
$\widetilde{\phi}(e^{i\theta})=\lim_{r\rightarrow1^-}\phi(re^{i\theta})$ a.e. on $\partial\D$,
and where it follows from the relation $g=\phi f_{_0}\in H^p(\D)$ that the radial limit $\lim_{r\rightarrow 1^-}f_{_0}(re^{i\theta})$ of
$f_{_0}(z)$ exists almost everywhere on $\partial\D$.

Now, we give a specific function $f_{_0}$ in $A^p_{\alpha}(\D)\setminus H^{^p}(\D)$ that
has a radial limit almost nowhere on $\partial\D$,
which contradicts the conclusion of the preceding paragraph.

In fact, let
$$f_{_0}(z)=\sum^{\infty}_{k=1}z^{3^{^{k-1}}}=z+z^3+z^9+z^{27}+\cdots+ z^{3^{^{k-1}}}+\cdots,$$
for all $z\in \D.$
Write $a_{_k}=1,~n_{_k}=3^{k-1},~~k=1,2,\cdots.$ Then it is
easy to see that
$$f_{_0}(z)=\sum^{\infty}_{k=1}a_{_k}z^{n_{_k}},~~\frac{n_{_{k+1}}}{n_{_k}}=3,~~\sum^{\infty}_{k=1}\frac{|a_{_k}|}{n^{\alpha+1}_{k}}< \infty,~~\mbox{and} ~~
\sum^{\infty}_{k=1}|a_{_k}|^2=\infty.$$
Thus by the criterion of $A^p_{\alpha}(\D)$-functions (see for instance \cite{rzz08}, $\S$ 14 Lacunary series),
it can be obtain that
$f_{_0} \in A^p_{\alpha}(\D)$ for any $1\leq p <\infty$ and any $-1<\alpha < \infty$,
and by the criterion of the radial limit (see for instance \cite{pda04} p.80),
the function $f_{_0}$ has a radial limit
almost nowhere on $\partial\D$. On the other hand, since each function $f$ in $H^p(\D)$ has the
radial limit almost everywhere on $\partial\D$, it follows that
$f_{_0}\in A_\alpha^{^p}(\D)\backslash H^{^p}(\D)$.
This completes the proof.

\msk
\noindent
\textbf{
Note 1.} It follows from the proof of Theorem 2 that for
any inner function $\phi$, $\phi H^{^p}(\D)$ is a linear subspace of the Bergman space
$A_\alpha^{^p}(\D)~(1\leq p<\infty, -1<\alpha <\infty)$, but is not a closed linear subspace in $A^{^p}_\alpha(\D)$ (in the norm of $A^{^p}_\alpha(\D)$).
In particular, $H^{^p}(\D)$ is a linear subspace of $A_\alpha^{^p}(\D)$, but it is not a closed linear subspace of $A_\alpha^{^p}(\D)$
(in the norm of $A^{^p}_\alpha(\D)$).

\msk

By the way, although every invariant subspace $M(\neq \{0\})$ for
the shift operator $M_z$ on the Hardy space $H^{^p}(\D)\ ~~(1\leq p<\infty)$ has the form
$M=\phi H^{p}(\D)$
for some inner function $\phi$,
and many invariant subspaces $M$ for
the shift operator $M_{e^{i\theta}}$ on the Lebesgue space $L^2(\partial\D)$ also has the form
$M=\widetilde{\phi} H^{2}(\partial\D)$
for some inner function $\widetilde{\phi}$,
but by Theorem 2, any invariant subspace $M$ for the shift operator $M_z$ on the Bergman space
$A_\alpha^{^p}(\D)~(1\leq p<\infty, -1<\alpha <\infty)$ has no the form $M=\phi H^{^p}(\D)$
for any inner function $\phi$.

\msk
\noindent
\textbf{
Note 2.} It can be seen from the proof of Theorem 2 that for each function $\phi$
in $H^\infty(\DD)$ (of course, for each inner function $\phi$), we have
\begin{equation}\label{002}
  \phi H^p(\D)\subsetneqq \phi A^p_\alpha{\D}
\end{equation}
for all $1\leq p<\infty$ and all $-1<\alpha <\infty$. In fact, it is well known that
$H^{^p}(\D) \subset A^{^p}_\alpha(\D)~(1\leq p<\infty, -1<\alpha <\infty)$. Therefore
$\phi H^{^p}(\D) \subset \phi A^{^p}_\alpha(\D)$.
Thus if the expression (\ref{002}) were not
true, then there would be a result that
\begin{equation}\label{003}
\phi H^p(\D)=\phi A^p_\alpha(\D)~{\mbox{(as two set)}}
\end{equation}
for some value of $p$ and some value of $\alpha$.
On the other hand, it has be shown in the
proof of Theorem 2 that the function
$$f_{_0}(z)=\sum^{\infty}_{k=1}z^{3^{^{k-1}}}=z+z^3+z^9+z^{27}+\cdots+ z^{3^{^{k-1}}}+\cdots$$
belongs to $A^p_{\alpha}(\D)$ for any $1\leq p<\infty$ and any $-1<\alpha <\infty$, and $f_{_0}$
has a radial limit almost nowhere on $\partial \D$.
Thus by (\ref{003}), there is a function $h_{_0}\in H^p(\D)$ such that $\phi h_{_0}=\phi f_{_0}$. Let
$g_{_0}=\phi h_{_0}$, then $g_{_0}\in H^p(\D)$ and $\phi f_{_0}=g_{_0}$. Therefore
$$\lim_{r\rightarrow 1^-}f_{_0}(re^{i\theta})=\lim_{r\rightarrow 1^-}
\frac{g_{_0}(re^{i\theta})}{\phi(re^{i\theta})}
=\frac{\widetilde{g_{_0}}(e^{i\theta})}{\widetilde{\phi}(e^{i\theta})}~~~~~~~ {\mbox{a.~e. ~on} ~~~ \partial\D},$$
where $\widetilde{g_{_0}}(e^{i\theta})=\lim_{r\rightarrow 1^-} g_{_0}(re^{i\theta})$
a. e. on $\partial\D$ and $\widetilde{\phi}(e^{i\theta})=\lim_{r\rightarrow 1^-} \phi(re^{i\theta})$
a. e. on $\partial\D$.
This contradicts the fact that $f_{_0}$ has a radial limit almost nowhere on $\partial \D$.

\end{document}